\ifx\documentclass\undefined
\documentstyle[12pt]{article}
\else
\documentclass[12pt]{article}
\usepackage{latexsym}
\usepackage{amsfonts}
\usepackage{amsmath}
\usepackage{amssymb}
\usepackage{graphicx}
\fi

\author{ K\'aroly Bezdek  
\thanks{Partially supported by a Natural Sciences and 
Engineering Research Council of Canada Discovery Grant.} }

\font\tenBbb=msbm10 at 12pt         \font\sevenBbb=msbm9    \font\fiveBbb=msbm7
\newfam\Bbbfam
\textfont\Bbbfam=\tenBbb \scriptfont\Bbbfam=\sevenBbb
\scriptscriptfont\Bbbfam=\fiveBbb

\def\kkk{\null\hfill $\Box$\smallskip}

\newtheorem{theorem}{Theorem}[section]
\newtheorem{lemma}[theorem]{Lemma}

\newtheorem{sled}[theorem]{Corollary}
\newtheorem{con}[theorem]{Conjecture}

\newcommand{\proof}{{\noindent\bf Proof:{\ \ }}}

\title{Contact numbers for congruent sphere packings in Euclidean $3$-space 
}

\begin{document}

\maketitle

\date

\begin{abstract}
Continuing the investigations of Harborth (1974)  and the author (2002) we study the 
following two rather basic problems on sphere packings. Recall that the contact graph of an arbitrary finite packing 
of unit balls (i.e., of an arbitrary finite family of non-overlapping unit balls) in Euclidean 
3-space is the (simple) graph whose vertices correspond to the packing elements and whose 
two vertices are connected by an edge if the corresponding two packing 
elements touch each other. One of the most basic questions on contact graphs is 
to find the maximum number of edges that a contact graph of a packing of n unit balls 
can have in Euclidean 3-space. Our method for finding lower and upper estimates for 
the largest contact numbers is a combination of analytic and combinatorial ideas and 
it is also based on some recent results on sphere packings. 
Finally, we are interested also in the following more special version of the above problem. 
Namely, let us imagine that we are given a lattice unit sphere packing with the center points 
forming the lattice $\Lambda$ in Euclidean 3-space (and with certain pairs of unit balls touching each other) and then let us generate packings of n unit balls such that each and every 
center of the n unit balls is chosen from $\Lambda$. Just as in the general case we are interested in finding good estimates for the largest contact number of the packings of n unit balls obtained in this way.

\end{abstract}

\section{Introduction}

Let $\mathbb{E}^{d}$ denote the $d$-dimensional Euclidean space. Then the {\it contact graph} of an arbitrary finite packing of unit balls (i.e., of an arbitrary finite family of non-overlapping unit balls) in $\mathbb{E}^{d}$ is the (simple) graph whose vertices correspond to the packing elements and whose two vertices are connected by an edge if and only if the corresponding two packing elements touch each other. One of the most basic questions on contact graphs is to find the maximum number of edges that a contact graph of a packing of $n$ unit balls can have in $\mathbb{E}^{d}$. In 1974 Harborth \cite{Ha} proved the following optimal result in $\mathbb{E}^{2}$: the maximum number $c(n)$ of touching pairs in a packing of $n$ congruent circular disks in $\mathbb{E}^{2}$ is precisely $\lfloor 3n-\sqrt{12n-3}\rfloor$ implying that $$\lim_{n\to +\infty}\frac{3n-c(n)}{\sqrt{n}}=\sqrt{12}=3.464\dots \ .$$ Some years later the author \cite{B02} has proved the following estimates in higher dimensions. The number of touching pairs in an arbitrary packing of $n>1$ unit balls in $\mathbb{E}^{d}$, $d\geq 3$ is less than
$$\frac{1}{2}\tau_d\, n-\frac{1}{2^{d}}\delta_d^{-\frac{d-1}{d}}\; n^{\frac{d-1}{d}},$$
where $\tau_d$ stands for the kissing number of a unit ball in $\mathbb{E}^{d}$ (i.e., it denotes the maximum number of non-overlapping unit balls of $\mathbb{E}^{d}$ that can touch a given unit ball in $\mathbb{E}^{d}$) and $\delta_d$ denotes the largest possible density for (infinite) packings of unit balls in $\mathbb{E}^{d}$. Now, recall that on the one hand, according to the well-known theorem of Kabatiansky and Levenshtein \cite{KL} $\tau_d\leq 2^{0.401d(1+\text{o}(1))}$ and $\delta_d\leq 2^{-0.599d(1+\text{o}(1))}$ as $d\to +\infty$ on the other hand, $\tau_3=12$ (for the first complete proof see \cite{SW}) moreover, according to the recent breakthrough result of Hales \cite{H05} $\delta_3=\frac{\pi}{\sqrt{18}}$. Thus, by combining the above results together we get that the number of touching pairs in an arbitrary packing of $n>1$ unit balls in $\mathbb{E}^{d}$ is less than $$\frac{1}{2}2^{0.401d(1+\text{o}(1))}\, n-\frac{1}{2}2^{-0.401(d-1)(1-\text{o}(1))}\; n^{\frac{d-1}{d}}$$ as $d\to +\infty$ and in particular, it is less than $$6n-\frac{1}{8}\left(\frac{\pi}{\sqrt{18}}\right)^{-\frac{2}{3}}\; n^{\frac{2}{3}}=6n-0.152\dots n^{\frac{2}{3}}$$ for $d=3$. The main purpose of this note is to improve further the latter result. In order, to state our theorem in a proper form we need to introduce a bit of additional terminology. If $\cal P$ is a packing of $n$ unit balls in $\mathbb{E}^{3}$, then let $C(\cal P)$ stand for the number of touching pairs in $\cal P$, that is, let $C(\cal P)$ denote the number of edges of the contact graph of $\cal P$ and call it the {\it contact number} of $\cal P$. Moreover, let $C(n)$ be the largest $C(\cal P)$ for packings $\cal P$ of $n$ unit balls in $\mathbb{E}^{3}$. Finally, let us imagine that we generate packings of n unit balls in $\mathbb{E}^{3}$ in such a special way that each and every center of the $n$ unit balls chosen, is a lattice point of some fixed lattice $\Lambda$ (resp., of the face-centered cubic lattice $\Lambda_{fcc}$) with shortest non-zero lattice vector of length $2$. (Here, a lattice means a (discrete) set of points having position vectors that are integer linear combinations of three fixed linearly independent vectors of $\mathbb{E}^{3}$.) Then let $C_{\Lambda}(n)$ (resp., $C_{fcc}(n)$) denote the largest possible contact number of all packings of $n$ unit balls obtained in this way. Before stating our main theorem we make the following comments. First, recall that according to \cite{H05} the lattice unit sphere packing generated by $\Lambda_{fcc}$ gives the largest possible density for unit ball packings in $\mathbb{E}^{3}$, namely $\frac{\pi}{\sqrt{18}}$ with each ball touched by $12$ others such that their centers form the vertices of a cuboctahedron. Second, it is easy to see that $C_{fcc}(2)=C(2)=1, C_{fcc}(3)=C(3)=3, C_{fcc}(4)=C(4)=6$. Third, it is natural to conjecture that $C_{fcc}(9)=C(9)=21$. Based on the trivial inequalities $C(n+1)\ge C(n)+3, C_{fcc}(n+1)\ge C_{fcc}(n)+3$ valid for all $n\ge 2$, it would follow that $C_{fcc}(5)=C(5)=9,C_{fcc}(6)=C(6)=12, C_{fcc}(7)=C(7)=15$, and $C_{fcc}(8)=C(8)=18$. Furthermore, we note that $C(10)\geq 25, C(11)\geq 29$, and $C(12)\geq 33$.  In order, to see that one should take the union $\mathbf{U}$ of two regular octahedra of edge length $2$ in $\mathbb{E}^{3}$ such that they share a regular triangle face $T$ in common and lie on opposite sides of it. If we take the unit balls centered at the nine vertices of $\mathbf{U}$, then there are exactly $21$ touching pairs among them. Also, we note that along each side of $T$ the dihedral angle of $\mathbf{U}$ is concave and in fact, it can be completed to $2\pi$ by adding twice the dihedral angle of a regular tetrahedron in $\mathbb{E}^{3}$. This means that along each side of $T$ two triangular faces of  $\mathbf{U}$ meet such that for their four vertices there exists precisely one point in $\mathbb{E}^{3}$ lying outside  $\mathbf{U}$ and at distance $2$ from each of the four vertices. Finally, if we take the twelve vertices of a cuboctahedron of edge length $2$ in $\mathbb{E}^{3}$ along with its center of symmetry, then the thirteen unit balls centered about them have $36$ contacts implying that $C(13)\geq 36$. Whether in any of the inequalities $C(10)\geq 25, C(11)\geq 29, C(12)\geq 33$, and $C(13)\geq 36$ we have equality is a challenging open question. In the rest of this note we give a proof of the following theorem.

\begin{theorem}\label{Bezdek-on-contact-numbers}

\item(i) \hskip0.6cm $C(n) < 6n-0.695n^{\frac{2}{3}}$
for all $n\ge 2$.

\item(ii) \hskip0.5cm Let $\Lambda$ be any lattice  of $\mathbb{E}^{3}$ with shortest non-zero lattice vector of length $2$. Then  $C_{\Lambda}(n)< 6n-\frac{3\sqrt[3]{18\pi}}{\pi}n^{\frac{2}{3}}=6n-3.665\dots n^{\frac{2}{3}}$
for all $n\ge 2$.

\item(iii) \hskip0.3cm $6n-\sqrt[3]{486}n^{\frac{2}{3}}<C_{fcc}(n)\leq C(n)$ for all $n=\frac{k(2k^2+1)}{3}$ with $k\ge 2$.

\textbf{}\end{theorem}

As an immediate result we get 

\begin{sled}
$$0.695<\frac{6n-C(n)}{n^{\frac{2}{3}}}< \sqrt[3]{486}=7.862\dots$$ 
for all $n=\frac{k(2k^2+1)}{3}$ with $k\ge 2$.
\end{sled}

The following was noted in \cite{B02}. Due to the Minkowski difference body method (see for example, Chapter 6 in \cite{Ro64}) the family ${\cal P}_{\mathbf{K}}:=\{\mathbf{t}_1+\mathbf{K}, \mathbf{t}_2+\mathbf{K}, \dots , \mathbf{t}_n+\mathbf{K}\} $ of $n$ translates of the convex body $\mathbf{K}$ in $\mathbb{E}^d$ is a packing if and only if the family ${\cal P}_{\mathbf{K}_{\mathbf{o}}}:=\{\mathbf{t}_1+\mathbf{K}_{\mathbf{o}}, \mathbf{t}_2+\mathbf{K}_{\mathbf{o}}, \dots , \mathbf{t}_n+\mathbf{K}_{\mathbf{o}}\} $ of $n$ translates of the symmetric difference body $\mathbf{K}_{\mathbf{o}}:=\frac{1}{2}(\mathbf{K}+(-\mathbf{K}))$ of $\mathbf{K}$ is a packing in $\mathbb{E}^d$. Moreover, the number of touching pairs in the packing ${\cal P}_{\mathbf{K}}$ is equal to the number of touching pairs in the packing ${\cal P}_{\mathbf{K}_{\mathbf{o}}}$. Thus, for this reason and for the reason that if $\mathbf{K}$ is a convex body of constant width in $\mathbb{E}^d$, then $\mathbf{K}_{\mathbf{o}}$ is a ball of $\mathbb{E}^d$, Theorem~\ref{Bezdek-on-contact-numbers} extends in a straightforward way to translative packings of convex bodies of constant width in $\mathbb{E}^3$.

For the sake of completeness we mention that the nature of contact numbers changes dramatically for non-congruent sphere packings in $\mathbb{E}^3$. For more details on that we refer the interested reader to the elegant paper \cite{KuOd94} of Kuperberg and Schramm.

Last but not least, it would be interesting to improve further the estimates of Theorem~\ref{Bezdek-on-contact-numbers}. In the last section of this paper we mention some ideas that (might) lead to such improvements for $(i)$ in Theorem~\ref{Bezdek-on-contact-numbers}.

\section{Proof of Theorem~\ref{Bezdek-on-contact-numbers}}
\bigskip

\subsection{Proof of (i)}

Let $\mathbf{B}$ denote the (closed) unit ball centered at the origin $\mathbf{o}$ of $\mathbb{E}^3$ and let ${\cal P}:=\{\mathbf{c}_1+\mathbf{B}, \mathbf{c}_2+\mathbf{B}, \dots , \mathbf{c}_n+\mathbf{B}\} $ denote the packing of $n$ unit balls with centers $\mathbf{c}_1, \mathbf{c}_2, \dots , \mathbf{c}_n$ in $\mathbb{E}^3$ having the largest number $C(n)$ of touching pairs among all packings of $n$ unit balls in $\mathbb{E}^3$. (${\cal P}$ might not be uniquely determined up to congruence in which case ${\cal P}$ stands for any of those extremal packings.) Now, let $\hat{r}:=1.81383$. The following statement shows the main property of $\hat{r}$ that is needed for our proof of Theorem~\ref{Bezdek-on-contact-numbers}.

\begin{lemma}\label{Bezdek-Boroczky-Szabo}
Let $\mathbf{B}_1, \mathbf{B}_2, \dots , \mathbf{B}_{13}$ be $13$ different members of a packing of unit balls in $\mathbb{E}^3$. Assume that each ball of the family $\mathbf{B}_2, \mathbf{B}_3, \dots , \mathbf{B}_{13}$ touches $\mathbf{B}_1$. Let $\hat{\mathbf{B}}_i$ be the closed ball concentric with $\mathbf{B}_i$ having radius $\hat{r}$, $1\le i\le 13$. Then
the boundary ${\rm bd}(\hat{\mathbf{B}}_1)$ of $\hat{\mathbf{B}}_1$ is covered by the balls $\hat{\mathbf{B}}_2, \hat{\mathbf{B}}_3, \dots , \hat{\mathbf{B}}_{13}$, that is,
$${\rm bd}(\hat{\mathbf{B}}_1)\subset \cup_{j=2}^{13}\hat{\mathbf{B}}_j \ .$$
\end{lemma}

\proof
The statement follows from the following recent result of B\"or\"oczky and Szab\'o
\cite{Bo2}.

\begin{lemma}\label{Boroczky-Szabo}
Let $\mathbf{B}_1, \mathbf{B}_2, \dots , \mathbf{B}_{14}$ be $14$ different members of a packing of unit balls in $\mathbb{E}^3$. Assume that each ball of the family $\mathbf{B}_2, \mathbf{B}_3, \dots , \mathbf{B}_{13}$ touches $\mathbf{B}_1$. Then the distance between
the centers of $\mathbf{B}_1$ and $\mathbf{B}_{14}$ is at least $$2.205279217705.$$
\end{lemma}

Indeed, Lemma~\ref{Boroczky-Szabo} combined with the following elementary trigonometry completes the proof of Lemma~\ref{Bezdek-Boroczky-Szabo}. The more exact details are as follows. Let $\alpha$ denote the measure of the angles opposite to the equal sides of the isosceles triangle $\triangle\mathbf{o}_1\mathbf{p}\mathbf{q}$ with ${\rm dist}(\mathbf{o}_1,\mathbf{p})=2$ and ${\rm dist}(\mathbf{p},\mathbf{q})={\rm dist}(\mathbf{o}_1,\mathbf{q})=\hat{r}$, where ${\rm dist}(\cdot , \cdot )$ denotes the Euclidean distance between the corresponding two points. Clearly, 
$\cos\alpha =\frac{1}{\hat{r}}$. Moreover, if $\triangle\mathbf{o}_1\mathbf{p}\mathbf{q'}$ denotes the isosceles triangle whose side $\mathbf{o}_1\mathbf{q'}$ contains $\mathbf{q}$ as a relative interior point such that ${\rm dist}(\mathbf{o}_1,\mathbf{p})={\rm dist}(\mathbf{p},\mathbf{q'})=2$, then ${\rm dist}(\mathbf{o}_1,\mathbf{q'})=4\cos\alpha=\frac{4}{\hat{r}}=2.205278333691... <2.205279217705$. This inequality together with Lemma~\ref{Boroczky-Szabo} implies in a straightforward way that if $\mathbf{o}_1$ is in fact, the center of the unit ball $\mathbf{B}_1$ touched by each of the unit balls $\mathbf{B}_2, \mathbf{B}_3, \dots , \mathbf{B}_{13}$ at the points $\mathbf{t}_j\in{\rm bd}(\mathbf{B}_j)\cap{\rm bd}(\mathbf{B}_1), 2\le j\le 13$, then the radius of the circumscribed circle of each face of the convex polyhedron ${\rm conv}\{\mathbf{t}_2, \mathbf{t}_3, \dots , \mathbf{t}_{13}\}$ is less than $\sin\alpha$ and so, ${\rm bd}(\hat{\mathbf{B}}_1)\subset \cup_{j=2}^{13}\hat{\mathbf{B}}_j$. This finishes the proof of Lemma~\ref{Bezdek-Boroczky-Szabo}. \kkk

Now, let us take the union $\bigcup_{i=1}^n\left(\mathbf{c}_i+\hat{r}\mathbf{B}\right)$ of the closed balls $\mathbf{c}_1+\hat{r}\mathbf{B}, \mathbf{c}_2+\hat{r}\mathbf{B}, \dots , \mathbf{c}_n+\hat{r}\mathbf{B}$ of radii $\hat{r}$ centered at the points $\mathbf{c}_1, \mathbf{c}_2, \dots , \mathbf{c}_n$ in $\mathbb{E}^3$. 

\begin{lemma}\label{Bezdek}
$$\frac{n{\rm vol}_3(\mathbf{B})}{{\rm vol}_3\left(\bigcup_{i=1}^n\left(\mathbf{c}_i+\hat{r}\mathbf{B}\right)\right)}<0.7785,$$
where ${\rm vol}_3(\cdot)$ refers to the $3$-dimensional volume of the corresponding set.
\end{lemma}

\proof
First, partition $\bigcup_{i=1}^n\left(\mathbf{c}_i+\hat{r}\mathbf{B}\right)$ into truncated Voronoi cells as follows. Let $\mathbf{P}_i$ denote the Voronoi cell of the packing $\cal P$ assigned to $\mathbf{c}_i+\mathbf{B}$, $1\le i\le n$, that is, let $\mathbf{P}_i$ stand for the set of points of $\mathbb{E}^{3}$ that are not farther away from $\mathbf{c}_i$ than from any other $\mathbf{c}_j$ with $j\neq i, 1\le j\le n$. Then, recall the well-known fact (see for example, \cite{Ro64}) that the Voronoi cells $\mathbf{P}_i$, $1\le i\le n$ just introduced form a tiling of $\mathbb{E}^{3}$. Based on this it is easy to see that the truncated Voronoi cells $\mathbf{P}_i\cap (\mathbf{c}_i+\hat{r}\mathbf{B})$, $1\le i\le n$ generate a tiling of the non-convex container $\bigcup_{i=1}^n\left(\mathbf{c}_i+\hat{r}\mathbf{B}\right)$ for the packing $\cal P$. Second, as $\sqrt{\frac{3}{2}}=1.2247...<\hat{r}=1.81383$ therefore the following recent result (Corollary 3 in \cite{B00}) of the author applied to the truncated Voronoi cells $\mathbf{P}_i\cap (\mathbf{c}_i+\hat{r}\mathbf{B})$, $1\le i\le n$ implies the inequality of Lemma~\ref{Bezdek} in a straightforward way.

\begin{lemma}\label{Bezdek-Voronoi}
Let $\cal F$ be an arbitrary (finite or infinite) family of non-overlap\-ping unit balls in $\mathbb{E}^3$ with the unit ball $\mathbf{B}$ centered at the origin $\mathbf{o}$ of $\mathbb{E}^3$ belonging to $\cal F$. Let $\mathbf{P}$ stand for the Voronoi cell of the packing $\cal F$ assigned to $\mathbf{B}$. Let $\mathbf{Q}$ denote a regular dodecahedron circumscribed $\mathbf{B}$ having circumradius $\sqrt{3}\tan\frac{\pi}{5}=1.2584...$. If $F$ denotes a pentagonal face of $\mathbf{Q}$, then let $F'$ be the regular pentagon positively homothetic to $F$ with respect to the center of $F$ lying in the plane of $F$ such that its vertices are at distance $\sqrt{1+\frac{1}{3\cos^2\frac{\pi}{5}}}=1.2285...$ from $\mathbf{o}$ (and its sides are at distance $\frac{2}{\sqrt{3}}=1.1547...$ from $\mathbf{o}$). Finally, let $r:=\sqrt{\frac{3}{2}}=1.2247...$ and $\mathbf{Q}':={\rm conv}\left((F'\cap r\mathbf{B})\cup\{\mathbf{o}\}\right)$, where $r\mathbf{B}$ denotes the ball of radius $r$ centered at the origin $\mathbf{o}$ of $\mathbb{E}^3$. Then
$$\frac{{\rm vol}_3(\mathbf{B})}{{\rm vol}_3(\mathbf{P})}\le\frac{{\rm vol}_3(\mathbf{B})}{{\rm vol}_3(\mathbf{P}\cap r\mathbf{B})}\le\frac{{\rm vol}_3(\mathbf{Q}'\cap \mathbf{B})}{{\rm vol}_3(\mathbf{Q}')}$$
$$= \frac{20\sqrt{6}\arctan(\frac{\sqrt{2}}{2})-2(2\sqrt{6}-1)\pi}{5\sqrt{2}+3\pi-15\arctan(\frac{\sqrt{2}}{2})}=0.77842...<0.7785 .$$

\end{lemma}

This finishes the proof of Lemma~\ref{Bezdek}. \kkk

The well-known isoperimetric inequality \cite{Os78} applied to $\bigcup_{i=1}^n\left(\mathbf{c}_i+\hat{r}\mathbf{B}\right)$ yields 

\begin{lemma}\label{isoperimetric-inequality}
$$36\pi{\rm vol}_3^2\left(\bigcup_{i=1}^n\left(\mathbf{c}_i+\hat{r}\mathbf{B}\right)\right)\le{\rm svol}_2^3\left({\rm bd}\left(\bigcup_{i=1}^n\left(\mathbf{c}_i+\hat{r}\mathbf{B}\right)\right)\right),$$
where ${\rm svol}_2(\cdot)$ refers to the $2$-dimensional surface volume of the corresponding set.
\end{lemma}

Thus, Lemma~\ref{Bezdek} and Lemma~\ref{isoperimetric-inequality} generate the following inequality.

\begin{sled}\label{lower-bound-for-surface-area-in-3D}
$$14.849236n^{\frac{2}{3}}<14.84923634...n^{\frac{2}{3}}=\frac{4\pi}{(0.7785)^{\frac{2}{3}}}n^{\frac{2}{3}}$$ $$< {\rm svol}_2\left({\rm bd}\left(\bigcup_{i=1}^n\left(\mathbf{c}_i+\hat{r}\mathbf{B}\right)\right)\right).$$
\end{sled}

Now, assume that $\mathbf{c}_i+\mathbf{B}\in {\cal P}$ is tangent to $\mathbf{c}_j+\mathbf{B}\in {\cal P}$ for all $j\in T_i$, where $T_i\subset\{1, 2, \dots , n\}$ stands for the family of indices $1\le j\le n$ for which ${\rm dist}(\mathbf{c}_i, \mathbf{c}_j)=2$. Then let
$\hat{S}_i:={\rm bd}(\mathbf{c}_i+\hat{r}\mathbf{B})$ and let $\hat{\mathbf{c}}_{ij}$ be the intersection of the line segment $\mathbf{c}_i\mathbf{c}_j$ with $\hat{S}_i$ for all $j\in T_i$. Moreover, let $C_{\hat{S}_i}(\hat{\mathbf{c}}_{ij}, \frac{\pi}{6})$ (resp., $C_{\hat{S}_i}(\hat{\mathbf{c}}_{ij}, \alpha)$) denote the open spherical cap of $\hat{S}_i$ centered at $\hat{\mathbf{c}}_{ij}\in \hat{S}_i$ having angular radius $\frac{\pi}{6}$ (resp., $\alpha$ with $0<\alpha<\frac{\pi}{2}$ and $\cos\alpha=\frac{1}{\hat{r}}$). Clearly, the family $\{C_{\hat{S}_i}(\hat{\mathbf{c}}_{ij}, \frac{\pi}{6}), j\in T_i\}$ consists of pairwise disjoint open spherical caps of $\hat{S}_i$; moreover,
\begin{equation}\label{Bezdek-estimate-I}
\frac{\sum_{j\in T_i}{\rm svol}_2\left(C_{\hat{S}_i}(\hat{\mathbf{c}}_{ij}, \frac{\pi}{6})\right)}{{\rm svol}_2\left(\cup_{j\in T_i}C_{\hat{S}_i}(\hat{\mathbf{c}}_{ij}, \alpha)\right)}=
\frac{\sum_{j\in T_i}{\rm Sarea}\left(C(\mathbf{u}_{ij}, \frac{\pi}{6})\right)}{{\rm Sarea}\left(\cup_{j\in T_i}C(\mathbf{u}_{ij}, \alpha)\right)},
\end{equation}
where $\mathbf{u}_{ij}:=\frac{1}{2}(\mathbf{c}_j-\mathbf{c}_i)\in \mathbb{S}^2:={\rm bd}(\mathbf{B})$ and $C(\mathbf{u}_{ij}, \frac{\pi}{6})\subset \mathbb{S}^2$ (resp., $C(\mathbf{u}_{ij}, \alpha)\subset \mathbb{S}^2$) denotes the open spherical cap of $\mathbb{S}^2$ centered at $\mathbf{u}_{ij}$ having angular radius $\frac{\pi}{6}$ (resp., $\alpha$)
and where ${\rm Sarea}(\cdot)$ refers to the spherical area measure on $\mathbb{S}^2$. Now, Moln\'ar's density bound (Satz I in \cite{Mo65}) implies that
\begin{equation}\label{Bezdek-estimate-II}
\frac{\sum_{j\in T_i}{\rm Sarea}\left(C(\mathbf{u}_{ij}, \frac{\pi}{6})\right)}{{\rm Sarea}\left(\cup_{j\in T_i}C(\mathbf{u}_{ij}, \alpha)\right)}<0.89332\ .
\end{equation}

In order to estimate ${\rm svol}_2\left({\rm bd}\left(\bigcup_{i=1}^n\left(\mathbf{c}_i+\hat{r}\mathbf{B}\right)\right)\right)$ from above let us assume that $m$ members of ${\cal P}$ have $12$ touching neighbours in ${\cal P}$ and $k$ members of ${\cal P}$ have at most $9$ touching neighbours in ${\cal P}$. Thus, $n-m-k$ members of ${\cal P}$ have either $10$ or $11$ touching neighbours in ${\cal P}$. (Here we have used the well-known fact that $\tau_3=12$, that is, no member of ${\cal P}$ can have more than $12$ touching neighbours.) Without loss of generality we may assume that $4\le k\le n-m$. 

First, we note that ${\rm Sarea}\left(C(\mathbf{u}_{ij}, \frac{\pi}{6})\right)=2\pi(1-\cos\frac{\pi}{6})=2\pi(1-\frac{\sqrt{3}}{2})$ and ${\rm svol}_2\left(C_{\hat{S}_i}(\hat{\mathbf{c}}_{ij}, \frac{\pi}{6})\right)=2\pi(1-\frac{\sqrt{3}}{2})\hat{r}^2$. Second, recall Lemma~\ref{Bezdek-Boroczky-Szabo} according to which if a member of ${\cal P}$ say, $\mathbf{c}_i+\mathbf{B}$ has exactly $12$ touching neighbours in ${\cal P}$, then $\hat{S}_i\subset \bigcup_{j\in T_i}(\mathbf{c}_j+\hat{r}\mathbf{B})$. These facts together with (\ref{Bezdek-estimate-I}) and (\ref{Bezdek-estimate-II}) imply the following estimate.

\begin{sled}\label{upper-bound-for-surface-area-in-3D}
$${\rm svol}_2\left({\rm bd}\left(\bigcup_{i=1}^n\left(\mathbf{c}_i+\hat{r}\mathbf{B}\right)\right)\right)
<\frac{32.04253}{3}(n-m-k)+32.04253k \ .
$$
\end{sled}

\proof
$${\rm svol}_2\left({\rm bd}\left(\bigcup_{i=1}^n\left(\mathbf{c}_i+\hat{r}\mathbf{B}\right)\right)\right)$$
$$<\left(4\pi\hat{r}^2-\frac{10\cdot2\pi(1-\frac{\sqrt{3}}{2})\hat{r}^2}{0.89332}\right)(n-m-k)+\left(4\pi\hat{r}^2-\frac{3\cdot2\pi(1-\frac{\sqrt{3}}{2})\hat{r}^2}{0.89332}\right)k$$
$$<10.34119(n-m-k)+32.04253k<\frac{32.04253}{3}(n-m-k)+32.04253k\ .$$ \kkk

Hence, Corollary~\ref{lower-bound-for-surface-area-in-3D} and Corollary~\ref{upper-bound-for-surface-area-in-3D} yield in a straightforward way that

\begin{equation}\label{Bezdek-estimate-III}
1.39026n^{\frac{2}{3}}-3k<n-m-k \ .
\end{equation}

Finally, as the number $C(n)$ of touching pairs in ${\cal P}$ is obviously at most $$
\frac{1}{2}\left( 12n-(n-m-k)-3k\right)\ ,
$$
therefore (\ref{Bezdek-estimate-III}) implies that
$$
C(n)\le \frac{1}{2}\left( 12n-(n-m-k)-3k\right) < 6n-0.69513n^{\frac{2}{3}}<6n-0.695n^{\frac{2}{3}},
$$
finishing the proof of $(i)$ in Theorem~\ref{Bezdek-on-contact-numbers}.

\subsection{Proof of (ii)}

First, we prove that  $C_{\Lambda}(n)\le C_{fcc}(n)$ and second, we give a proof of the desired upper bound for $C_{fcc}(n)$. The details are as follows.

Recall Voronoi's theorem (see \cite{CS92}) according to which every $3$-dimensional lattice is of the {\it first kind} i.e., it has an {\it obtuse superbase}. Thus, for the lattice $\Lambda$ (resp., $\Lambda_{fcc}$) we have a set of vectors $\mathbf{v}_0, \mathbf{v}_1, \mathbf{v}_2, \mathbf{v}_3$ (resp., $\mathbf{w}_0, \mathbf{w}_1, \mathbf{w}_2, \mathbf{w}_3$ ) such that $\mathbf{v}_1, \mathbf{v}_2, \mathbf{v}_3$ (resp., $\mathbf{w}_1, \mathbf{w}_2, \mathbf{w}_3$) is an integral basis for $\Lambda$ (resp., $\Lambda_{fcc}$) and $\mathbf{v}_0+\mathbf{v}_1+\mathbf{v}_2+\mathbf{v}_3=\mathbf{o}$ (resp., $\mathbf{w}_0+\mathbf{w}_1+\mathbf{w}_2+\mathbf{w}_3=\mathbf{o}$), and in addition $\mathbf{v}_i\cdot\mathbf{v}_j\le 0$ (resp., $\mathbf{w}_i\cdot\mathbf{w}_j\le 0$) for all $i,j=0,1,2,3$, $i\neq j$. Here ${}\cdot{}$ refers to the standard inner product of $\mathbb{E}^3$. Let $\mathbf{P}$ (resp., $\mathbf{Q}$) denote the Voronoi cell for the origin $\mathbf{o}\in \Lambda$ (resp., $\mathbf{o}\in \Lambda_{fcc}$) consisting of points of $\mathbb{E}^3$ that are at least as close to $\mathbf{o}$ as to any other lattice point of $\Lambda$ (resp., $\Lambda_{fcc}$). A vector $\mathbf{v}\in \Lambda$ (resp., $\mathbf{w}\in \Lambda_{fcc}$) is called a {\it strict Voronoi vector} of $\Lambda$ (resp., $\Lambda_{fcc}$) if the plane $\{\mathbf{x}\in \mathbb{E}^3\ | \ \mathbf{x}\cdot\mathbf{v}=\frac{1}{2}\mathbf{v}\cdot\mathbf{v}\} $ (resp., $\{\mathbf{x}\in \mathbb{E}^3\ | \ \mathbf{x}\cdot\mathbf{w}=\frac{1}{2}\mathbf{w}\cdot\mathbf{w}\} $ ) intersects $\mathbf{P}$ (resp., $\mathbf{Q}$) in a face. We need the following claim proved in \cite{CS92}. The list of $14$ lattice vectors of $\Lambda$ (resp., $\Lambda_{fcc}$) consisting of $$\pm\mathbf{v}_1,\pm(\mathbf{v}_0+\mathbf{v}_1), \pm(\mathbf{v}_1+\mathbf{v}_2), \pm(\mathbf{v}_1+\mathbf{v}_3),$$
$$\pm(\mathbf{v}_0+\mathbf{v}_1+\mathbf{v}_2), \pm(\mathbf{v}_0+\mathbf{v}_1+\mathbf{v}_3), \pm(\mathbf{v}_1+\mathbf{v}_2+\mathbf{v}_3)  $$
$${\rm (resp.,}\pm\mathbf{w}_1,\pm(\mathbf{w}_0+\mathbf{w}_1), \pm(\mathbf{w}_1+\mathbf{w}_2), \pm(\mathbf{w}_1+\mathbf{w}_3),$$
$$\pm(\mathbf{w}_0+\mathbf{w}_1+\mathbf{w}_2), \pm(\mathbf{w}_0+\mathbf{w}_1+\mathbf{w}_3), \pm(\mathbf{w}_1+\mathbf{w}_2+\mathbf{w}_3)  {\rm )}$$  
includes all the strict Voronoi vectors of $\Lambda$ (resp., $\Lambda_{fcc}$). As is well known (and in fact, it is easy check) at most 12 (resp., exactly 12) of the above 14 vectors has length $2$ and the others are of length strictly greater than $2$. Thus, it follows that without loss of generality we may assume that whenever $\mathbf{v}_i\cdot\mathbf{v}_i=4$ holds we have $\mathbf{w}_i\cdot\mathbf{w}_i=4$ as well. This implies the exisctence of a map $f : \Lambda\rightarrow\Lambda_{fcc}$ with the property that if ${\rm dist}(\mathbf{x},\mathbf{y})=2$ with $\mathbf{x}, \mathbf{y}\in \Lambda$, then also ${\rm dist}(f(\mathbf{x}),f(\mathbf{y}))=2$ holds. Indeed, $f$ can be defined via $f(\alpha\mathbf{v}_1+\beta\mathbf{v}_2+\gamma\mathbf{v}_3)=
\alpha\mathbf{w}_1+\beta\mathbf{w}_2+\gamma\mathbf{w}_3$ with $\alpha, \beta, \gamma $ being arbitrary integers. As a result we get the following: if ${\cal P}$ is a packing of $n$ unit balls with centers $\mathbf{c}_1, \mathbf{c}_2, \dots , \mathbf{c}_n\in \Lambda$, then the packing ${\cal P}_f$ of $n$ unit balls centered at the points $f(\mathbf{c}_1), f(\mathbf{c}_2), \dots , f(\mathbf{c}_n)\in \Lambda_{fcc}$ possesses the property that $C({\cal P})\le C({\cal P}_f)$. Thus, indeed, $C_{\Lambda}(n)\le C_{fcc}(n)$.

Although the idea of the proof of $(ii)$ for $C_{fcc}(n)$ is similar to that of $(i)$ they differ in the combinatorial counting part (see (\ref{Bezdek-estimate-X})) as well as in the density estimate for packings of spherical caps of angular radii $\frac{\pi}{6}$ (see (\ref{Bezdek-estimate-IX})). Moreover, the proof of $(ii)$ is based on the new parameter value $\bar{r}:=\sqrt{2}$ (replacing $\hat{r}=1.81383$). The details are as follows.

First, recall that if $\Lambda_{fcc}$ denotes the face-centered cubic lattice with shortest non-zero lattice vector of length $2$ in $\mathbb{E}^3$ and we place unit balls centered at each lattice point of $\Lambda_{fcc}$, then we get the fcc lattice packing of unit balls, labelled by ${\cal P}_{fcc}$, in which each unit ball is touched by $12$ others such that their centers form the vertices of a cuboctahedron. (Recall that a cuboctahedron is a convex polyhedron with $8$ triangular faces and $6$ square faces having $12$ identical vertices, with $2$ triangles and $2$ squares meeting at each, and $24$ identical edges, each separating a triangle from a square. As such it is a quasiregular polyhedron, i.e. an Archimedean solid, being vertex-transitive and edge-transitive.) Second, it is well-known (see \cite{F64} for more details) that the Voronoi cell of each unit ball in ${\cal P}_{fcc}$ is a rhombic dodecahedron (the dual of a cuboctahedron) of volume $\sqrt{32}$ and thus, the density of ${\cal P}_{fcc}$ is $\frac{\pi}{\sqrt{18}}$.

Now, let $\mathbf{B}$ denote the unit ball centered at the origin $\mathbf{o}\in\Lambda_{fcc}$ of $\mathbb{E}^3$ and let ${\cal P}:=\{\mathbf{c}_1+\mathbf{B}, \mathbf{c}_2+\mathbf{B}, \dots , \mathbf{c}_n+\mathbf{B}\} $ denote the packing of $n$ unit balls with centers $\{\mathbf{c}_1, \mathbf{c}_2, \dots , \mathbf{c}_n\}\subset\Lambda_{fcc}$ in $\mathbb{E}^3$ having the largest number $C_{fcc}(n)$ of touching pairs among all packings of $n$ unit balls being a sub-packing of ${\cal P}_{fcc}$. (${\cal P}$ might not be uniquely determined up to congruence in which case ${\cal P}$ stands for any of those extremal packings.)

The following two facts follow from the above description of ${\cal P}_{fcc}$ in a straightforward way. Let $\mathbf{B}_1, \mathbf{B}_2, \dots , \mathbf{B}_{13}$ be $13$ different members of ${\cal P}_{fcc}$ such that each ball of the family $\mathbf{B}_2, \mathbf{B}_3, \dots , \mathbf{B}_{13}$ touches $\mathbf{B}_1$. Moreover, let $\bar{\mathbf{B}}_i$ be the closed ball concentric with $\mathbf{B}_i$ having radius $\bar{r}=\sqrt{2}$, $1\le i\le 13$. Then
the boundary ${\rm bd}(\bar{\mathbf{B}}_1)$ of $\bar{\mathbf{B}}_1$ is covered by the balls $\bar{\mathbf{B}}_2, \bar{\mathbf{B}}_3, \dots , \bar{\mathbf{B}}_{13}$, that is,

\begin{equation}\label{Bezdek-estimate-IV}
{\rm bd}(\bar{\mathbf{B}}_1)\subset \cup_{j=2}^{13}\bar{\mathbf{B}}_j \ .
\end{equation}
In fact, $\bar{r}$ is the smallest radius with the above property. Moreover,

\begin{equation}\label{Bezdek-estimate-V}
\frac{n{\rm vol}_3(\mathbf{B})}{{\rm vol}_3(\bigcup_{i=1}^n\left(\mathbf{c}_i+\bar{r}\mathbf{B}\right))}<\frac{\pi}{\sqrt{18}}=0.7404...\ .
\end{equation}

As a next step we apply the isoperimetric inequality (\cite{Os78}): 

\begin{equation}\label{Bezdek-estimate-VI}
36\pi{\rm vol}_3^2\left(\bigcup_{i=1}^n\left(\mathbf{c}_i+\bar{r}\mathbf{B}\right)\right)\le{\rm svol}_2^3\left({\rm bd}\left(\bigcup_{i=1}^n\left(\mathbf{c}_i+\bar{r}\mathbf{B}\right)\right)\right) \ .
\end{equation}

Thus, (\ref{Bezdek-estimate-V}) and (\ref{Bezdek-estimate-VI}) yield in a straightforward way that

\begin{equation}\label{Bezdek-estimate-VII}
15.3532...n^{\frac{2}{3}}=4\sqrt[3]{18\pi}n^{\frac{2}{3}}
<{\rm svol}_2\left({\rm bd}\left(\bigcup_{i=1}^n\left(\mathbf{c}_i+\bar{r}\mathbf{B}\right)\right)\right) \ .
\end{equation}

Now, assume that $\mathbf{c}_i+\mathbf{B}\in {\cal P}$ is tangent to $\mathbf{c}_j+\mathbf{B}\in {\cal P}$ for all $j\in T_i$, where $T_i\subset\{1, 2, \dots , n\}$ stands for the family of indices $1\le j\le n$ for which ${\rm dist}(\mathbf{c}_i, \mathbf{c}_j)=2$. Then let
$\bar{S}_i:={\rm bd}(\mathbf{c}_i+\bar{r}\mathbf{B})$ and let $\bar{\mathbf{c}}_{ij}$ be the intersection of the line segment $\mathbf{c}_i\mathbf{c}_j$ with $\bar{S}_i$ for all $j\in T_i$. Moreover, let $C_{\bar{S}_i}(\bar{\mathbf{c}}_{ij}, \frac{\pi}{6})$ (resp., $C_{\bar{S}_i}(\bar{\mathbf{c}}_{ij}, \frac{\pi}{4})$) denote the open spherical cap of $\bar{S}_i$ centered at $\bar{\mathbf{c}}_{ij}\in \bar{S}_i$ having angular radius $\frac{\pi}{6}$ (resp., $\frac{\pi}{4}$). Clearly, the family $\{C_{\bar{S}_i}(\bar{\mathbf{c}}_{ij}, \frac{\pi}{6}), j\in T_i\}$ consists of pairwise disjoint open spherical caps of $\bar{S}_i$; moreover,
\begin{equation}\label{Bezdek-estimate-VIII}
\frac{\sum_{j\in T_i}{\rm svol}_2\left(C_{\bar{S}_i}(\bar{\mathbf{c}}_{ij}, \frac{\pi}{6})\right)}{{\rm svol}_2\left(\cup_{j\in T_i}C_{\bar{S}_i}(\bar{\mathbf{c}}_{ij}, \frac{\pi}{4})\right)}=
\frac{\sum_{j\in T_i}{\rm Sarea}\left(C(\mathbf{u}_{ij}, \frac{\pi}{6})\right)}{{\rm Sarea}\left(\cup_{j\in T_i}C(\mathbf{u}_{ij}, \frac{\pi}{4})\right)},
\end{equation}
where $\mathbf{u}_{ij}=\frac{1}{2}(\mathbf{c}_j-\mathbf{c}_i)\in \mathbb{S}^2$ and $C(\mathbf{u}_{ij}, \frac{\pi}{6})\subset \mathbb{S}^2$ (resp., $C(\mathbf{u}_{ij}, \frac{\pi}{4})\subset \mathbb{S}^2$) denotes the open spherical cap of $\mathbb{S}^2$ centered at $\mathbf{u}_{ij}$ having angular radius $\frac{\pi}{6}$ (resp., $\frac{\pi}{4}$). Now, the geometry of the cuboctahedron representing the $12$ touching neighbours of an arbitrary unit ball in ${\cal P}_{fcc}$ implies in a straightforward way that

\begin{equation}\label{Bezdek-estimate-IX}
\frac{\sum_{j\in T_i}{\rm Sarea}\left(C(\mathbf{u}_{ij}, \frac{\pi}{6})\right)}{{\rm Sarea}\left(\cup_{j\in T_i}C(\mathbf{u}_{ij}, \frac{\pi}{4})\right)}\le 6(1-\frac{\sqrt{3}}{2})=0.8038... 
\end{equation}
with equality when $12$ spherical caps of angular radius $\frac{\pi}{6}$ are packed on $\mathbb{S}^2$.

Finally, as ${\rm Sarea}\left(C(\mathbf{u}_{ij}, \frac{\pi}{6})\right)=2\pi(1-\cos\frac{\pi}{6})$ and ${\rm svol}_2\left(C_{\bar{S}_i}(\bar{\mathbf{c}}_{ij}, \frac{\pi}{6})\right)=2\pi(1-\frac{\sqrt{3}}{2})\bar{r}^2$ therefore (\ref{Bezdek-estimate-VIII}) and (\ref{Bezdek-estimate-IX}) yield that

$$
{\rm svol}_2\left({\rm bd}\left(\bigcup_{i=1}^n\mathbf{c}_i+\bar{r}\mathbf{B}\right)\right)
\le 4\pi\bar{r}^2 n-\frac{1}{6(1-\frac{\sqrt{3}}{2})}2\left(2\pi\left(1-\frac{\sqrt{3}}{2}\right)\bar{r}^2\right)C_{fcc}(n)
$$

\begin{equation}\label{Bezdek-estimate-X}
=8\pi n-\frac{4\pi}{3}C_{fcc}(n) \ .
\end{equation}

Thus, (\ref{Bezdek-estimate-VII}) and (\ref{Bezdek-estimate-X}) imply that

\begin{equation}\label{Bezdek-estimate-XI}
4\sqrt[3]{18\pi}n^{\frac{2}{3}}< 8\pi n-\frac{4\pi}{3}C_{fcc}(n) \ .
\end{equation}

From (\ref{Bezdek-estimate-XI}) the inequality $C_{fcc}(n) < 6n-\frac{3\sqrt[3]{18\pi}}{\pi}n^{\frac{2}{3}}=6n-3.665\dots n^{\frac{2}{3}}$ follows in a straightforward way for all $n\ge 2$. This completes the proof of $(ii)$ in Theorem~\ref{Bezdek-on-contact-numbers}.

\subsection{Proof of (iii)}

It is rather easy to show that for any positive integer $k\ge 2$ there are $n(k):=\frac{2k^3+k}{3}=\frac{k(2k^2+1)}{3}$ lattice points of the face-centered cubic lattice $\Lambda_{fcc}$ such that their convex hull is a regular octahedron $\mathbf{K}\subset\mathbb{E}^3$ of edge length $2(k-1)$ having exactly $k$ lattice points along each of its edges. Now, draw a unit ball around each lattice point of $\Lambda_{fcc}\cap\mathbf{K}$ and label the packing of the $n(k)$ unit balls obtained in this way by ${\cal P}_{fcc}(k)$. It is easy to check that if the center of a unit ball of ${\cal P}_{fcc}(k)$ is a relative interior point of an edge (resp., of a face) of $\mathbf{K}$, then the unit ball in question has $7$ (resp., $9$) touching neighbours in ${\cal P}_{fcc}(k)$. Last but not least, any unit ball of ${\cal P}_{fcc}(k)$ whose center is an interior pont of $\mathbf{K}$ has $12$ touching neighbours in ${\cal P}_{fcc}(k)$. Thus, the contact number $C\left({\cal P}_{fcc}(k)\right)$ of the packing ${\cal P}_{fcc}(k)$ is equal to

$$
6\frac{2(k-2)^3+(k-2)}{3}+
36\frac{(k-3)^2+(k-3)}{2}
+42(k-2)+12=4k^3-6k^2+2k. 
$$

As a result we get that

\begin{equation}\label{Bezdek-estimate-XII}
C\left({\cal P}_{fcc}(k)\right)=6n(k)-6k^2 .
\end{equation}

Finally, as $\frac{2k^3}{3}<n(k)$ therefore $6k^2<\sqrt[3]{486}n^{\frac{2}{3}}(k)$ and so,
(\ref{Bezdek-estimate-XII}) implies $(iii)$ of Theorem~\ref{Bezdek-on-contact-numbers} in a straightforward way. 

\section{Prospects for improvements}

\subsection{Sharpening the estimate $(i)$ in Theorem~\ref{Bezdek-on-contact-numbers}}

The following improvement on the estimate of Lemma~\ref{Boroczky-Szabo} has just been announced in \cite{H11} (see Lemma 9.15 on p. 229).

\begin{lemma}\label{Hales-2011}
Let $\mathbf{B}_1, \mathbf{B}_2, \dots , \mathbf{B}_{14}$ be $14$ different members of a packing of unit balls in $\mathbb{E}^3$. Assume that each ball of the family $\mathbf{B}_2, \mathbf{B}_3, \dots , \mathbf{B}_{13}$ touches $\mathbf{B}_1$. Then the distance between
the centers of $\mathbf{B}_1$ and $\mathbf{B}_{14}$ is at least $2.52$.
\end{lemma}

Naturally, this statement can be used in the same way as Lemma~\ref{Boroczky-Szabo} to get the following imrovement for Lemma~\ref{Bezdek-Boroczky-Szabo}.

\begin{lemma}\label{Bezdek-Hales}
Let $\mathbf{B}_1, \mathbf{B}_2, \dots , \mathbf{B}_{13}$ be $13$ different members of a packing of unit balls in $\mathbb{E}^3$. Assume that each ball of the family $\mathbf{B}_2, \mathbf{B}_3, \dots , \mathbf{B}_{13}$ touches $\mathbf{B}_1$. Let $\hat{\mathbf{B}}_i$ be the closed ball concentric with $\mathbf{B}_i$ having radius $\hat{r}:=1.58731$, $1\le i\le 13$. Then
the boundary ${\rm bd}(\hat{\mathbf{B}}_1)$ of $\hat{\mathbf{B}}_1$ is covered by the balls $\hat{\mathbf{B}}_2, \hat{\mathbf{B}}_3, \dots , \hat{\mathbf{B}}_{13}$, that is,
$${\rm bd}(\hat{\mathbf{B}}_1)\subset \cup_{j=2}^{13}\hat{\mathbf{B}}_j \ .$$
\end{lemma}

Also, it has has just been announced in \cite{H11} (see Lemma 9.13 on p. 228) that Lemma~\ref{Bezdek-Voronoi} can be improved as follows.

\begin{lemma}\label{Hales-Voronoi}
Let $\cal F$ be an arbitrary (finite or infinite) family of non-overlap\-ping unit balls in $\mathbb{E}^3$ with the unit ball $\mathbf{B}$ centered at the origin $\mathbf{o}$ of $\mathbb{E}^3$ belonging to $\cal F$. Let $\mathbf{P}$ stand for the Voronoi cell of the packing $\cal F$ assigned to $\mathbf{B}$. Let $\mathbf{Q}$ denote a regular dodecahedron circumscribed $\mathbf{B}$ (having circumradius $\sqrt{3}\tan\frac{\pi}{5}=1.2584...$).  Finally, let $r:=\sqrt{2}=1.4142...$ and let $r\mathbf{B}$ denote the ball of radius $r$ centered at the origin $\mathbf{o}$ of $\mathbb{E}^3$. Then
$$\frac{{\rm vol}_3(\mathbf{B})}{{\rm vol}_3(\mathbf{P})}\le\frac{{\rm vol}_3(\mathbf{B})}{{\rm vol}_3(\mathbf{P}\cap r\mathbf{B})}\le\frac{{\rm vol}_3(\mathbf{B})}{{\rm vol}_3(\mathbf{Q})}<0.7547.$$
\end{lemma}

Now, let ${\cal P}:=\{\mathbf{c}_1+\mathbf{B}, \mathbf{c}_2+\mathbf{B}, \dots , \mathbf{c}_n+\mathbf{B}\} $ denote the packing of $n$ unit balls with centers $\mathbf{c}_1, \mathbf{c}_2, \dots , \mathbf{c}_n$ in $\mathbb{E}^3$ having the largest number $C(n)$ of touching pairs among all packings of $n$ unit balls in $\mathbb{E}^3$. Moreover, let us take the union $\bigcup_{i=1}^n\left(\mathbf{c}_i+\hat{r}\mathbf{B}\right)$ of the closed balls $\mathbf{c}_1+\hat{r}\mathbf{B}, \mathbf{c}_2+\hat{r}\mathbf{B}, \dots , \mathbf{c}_n+\hat{r}\mathbf{B}$ of radii $\hat{r}=1.58731$ centered at the points $\mathbf{c}_1, \mathbf{c}_2, \dots , \mathbf{c}_n$ in $\mathbb{E}^3$. As
$\sqrt{2}<\hat{r}$ therefore following the proof of Lemma~\ref{Bezdek} based on Lemma~\ref{Bezdek-Voronoi}, Lemma~\ref{Hales-Voronoi} implies the following better estimate for Lemma~\ref{Bezdek}.

\begin{lemma}\label{Hales-II-2011}
$$\frac{n{\rm vol}_3(\mathbf{B})}{{\rm vol}_3\left(\bigcup_{i=1}^n\left(\mathbf{c}_i+\hat{r}\mathbf{B}\right)\right)}<0.7547.$$
\end{lemma}

Thus, Lemma~\ref{Hales-II-2011} combined with Lemma~\ref{isoperimetric-inequality} improves the estimate of Corollary~\ref{lower-bound-for-surface-area-in-3D} as follows.

\begin{sled}\label{Hales-lower-bound-for-surface-area-in-3D}
$$15.159805n^{\frac{2}{3}}< {\rm svol}_2\left({\rm bd}\left(\bigcup_{i=1}^n\left(\mathbf{c}_i+\hat{r}\mathbf{B}\right)\right)\right).$$
\end{sled}
Finally, following the method of the proof of Corollary~\ref{upper-bound-for-surface-area-in-3D} and replacing Lemma~\ref{Bezdek-Boroczky-Szabo} by the stronger Lemma~\ref{Bezdek-Hales} we get an improved version of Corollary~\ref{upper-bound-for-surface-area-in-3D}.

\begin{sled}\label{upper-bound-for-surface-area-in-3D-improved-via-Hales}
$${\rm svol}_2\left({\rm bd}\left(\bigcup_{i=1}^n\left(\mathbf{c}_i+\hat{r}\mathbf{B}\right)\right)\right)
<\frac{24.53902}{3}(n-m-k)+24.53902k \ .
$$
\end{sled}

Hence, Corollary~\ref{Hales-lower-bound-for-surface-area-in-3D} and Corollary~\ref{upper-bound-for-surface-area-in-3D-improved-via-Hales} yield in a straightforward way that

\begin{equation}\label{Bezdek-Hales-estimate-III}
1.85335n^{\frac{2}{3}}-3k<n-m-k \ .
\end{equation}

Finally, as the number $C(n)$ of touching pairs in ${\cal P}$ is obviously at most $$
\frac{1}{2}\left( 12n-(n-m-k)-3k\right)\ ,
$$
therefore (\ref{Bezdek-Hales-estimate-III}) implies that
$$
C(n)\le \frac{1}{2}\left( 12n-(n-m-k)-3k\right) < 6n-0.926675n^{\frac{2}{3}}<6n-0.926n^{\frac{2}{3}}.
$$

\subsection{Another approach for improving $(i)$ in Theorem~\ref{Bezdek-on-contact-numbers}}

Let $\delta(\mathbf{K})$ denote the largest density of packings of translates of the convex body $\mathbf{K}$ in $\mathbb{E}^{d}$, $d\geq 3$. The following result has been proved by the author in \cite{B02}.

\begin{lemma}\label{Bezdek-Wills}
Let $\mathbf{K_o}$ be a convex body in $\mathbb{E}^{d}$, $d\geq 2$ symmetric about the origin $\mathbf{o}$ of $\mathbb{E}^{d}$ and let $\{\mathbf{c}_1+\mathbf{K_o}, \mathbf{c}_2+\mathbf{K_o}, \dots , \mathbf{c}_n+\mathbf{K_o}\}$ be an arbitrary packing
of $n\ge1$ translates of $\mathbf{K_o}$ in $\mathbb{E}^{d}$. Then
$$ \frac{n{\rm vol}_d(\mathbf{K_o})}{{\rm vol}_d(\bigcup_{i=1}^n(\mathbf{c}_i+2\mathbf{K_o}))}< \delta(\mathbf{K_o}).$$
\end{lemma}

Let $\mathbf{B}$ denote the unit ball centered at the origin $\mathbf{o}$ of $\mathbb{E}^3$ and let ${\cal P}:=\{\mathbf{c}_1+\mathbf{B}, \mathbf{c}_2+\mathbf{B}, \dots , \mathbf{c}_n+\mathbf{B}\} $ denote the packing of $n$ unit balls with centers $\mathbf{c}_1, \mathbf{c}_2, \dots , \mathbf{c}_n$ having the largest number $C(n)$ of touching pairs among all packings of $n$ unit balls in $\mathbb{E}^3$. (${\cal P}$ might not be uniquely determined up to congruence in which case ${\cal P}$ stands for any of those extremal packings.) The well-known  result of Hales \cite{H05} according to which $\delta_3=\frac{\pi}{\sqrt{18}}$ and Lemma~\ref{Bezdek-Wills} imply in a straightforward way

\begin{lemma}\label{Bezdek-Hales-Ferguson}
$$\frac{n{\rm vol}_3(\mathbf{B})}{{\rm vol}_3(\bigcup_{i=1}^n(\mathbf{c}_i+2\mathbf{B}))}< \delta(\mathbf{B})=\frac{\pi}{\sqrt{18}}.$$
\end{lemma}

The isoperimetric inequality \cite{Os78} yields 

\begin{lemma}\label{Burago-Zalgaller}
$$36\pi{\rm vol}_3^2\left(\bigcup_{i=1}^n(\mathbf{c}_i+2\mathbf{B})\right)\le{\rm svol}_2^3\left({\rm bd}\left(\bigcup_{i=1}^n(\mathbf{c}_i+2\mathbf{B})\right)\right).$$
\end{lemma}

Thus, Lemma~\ref{Bezdek-Hales-Ferguson} and Lemma~\ref{Burago-Zalgaller} generate the following inequality.

\begin{sled}\label{third-lower-bound-for-surface-area-in-3D}
$$4\sqrt[3]{18\pi}n^{\frac{2}{3}}< {\rm svol}_2\left({\rm bd}\left(\bigcup_{i=1}^n(\mathbf{c}_i+2\mathbf{B})\right)\right).$$
\end{sled}

Now, assume that $\mathbf{c}_i+\mathbf{B}\in {\cal P}$ is tangent to $\mathbf{c}_j+\mathbf{B}\in {\cal P}$ for all $j\in T_i$, where $T_i\subset\{1, 2, \dots , n\}$ stands for the family of indices $1\le j\le n$ for which $\|\mathbf{c}_i-\mathbf{c}_j\|=2$. Then let
$S_i:={\rm bd}(\mathbf{c}_i+2\mathbf{B})$ and let $C_{S_i}(\mathbf{c}_j, \frac{\pi}{6})$ denote the open spherical cap of $S_i$ centered at $\mathbf{c}_j\in S_i$ having angular radius $\frac{\pi}{6}$. Clearly, the family $\{C_{S_i}(\mathbf{c}_j, \frac{\pi}{6}), j\in T_i\}$ consists of pairwise disjoint open spherical caps of $S_i$; moreover,
$$
\frac{\sum_{j\in T_i}{\rm svol}_2\left(C_{S_i}(\mathbf{c}_j, \frac{\pi}{6})\right)}{{\rm svol}_2\left(\cup_{j\in T_i}C_{S_i}(\mathbf{c}_j, \frac{\pi}{3})\right)}=
\frac{\sum_{j\in T_i}{\rm Sarea}\left(C(\mathbf{u}_{ij}, \frac{\pi}{6})\right)}{{\rm Sarea}\left(\cup_{j\in T_i}C(\mathbf{u}_{ij}, \frac{\pi}{3})\right)},
$$
where $\mathbf{u}_{ij}:=\frac{1}{2}(\mathbf{c}_j-\mathbf{c}_i)\in \mathbb{S}^2$ and $C(\mathbf{u}_{ij}, \frac{\pi}{6})\subset \mathbb{S}^2$ (resp., $C(\mathbf{u}_{ij}, \frac{\pi}{3})\subset \mathbb{S}^2$) denotes the open spherical cap of $\mathbb{S}^2$ centered at $\mathbf{u}_{ij}$ having angular radius $\frac{\pi}{6}$ (resp., $\frac{\pi}{3}$). Now, we are ready to state the main conjecture of this section.

\begin{con}\label{Bezdek-Molnar-lattice-estimate-II}
$$
\frac{\sum_{j\in T_i}{\rm Sarea}\left(C(\mathbf{u}_{ij}, \frac{\pi}{6})\right)}{{\rm Sarea}\left(\cup_{j\in T_i}C(\mathbf{u}_{ij}, \frac{\pi}{3})\right)}\le 6\left(1-\frac{\sqrt{3}}{2} \right)=0.8038\dots\ ,
$$
with equality when $12$ spherical caps of angular radius $\frac{\pi}{6}$ are packed on
$\mathbb{S}^2$.
\end{con}

If true, Conjecture~\ref{Bezdek-Molnar-lattice-estimate-II} can be used to improve further the upper bound for $C(n)$ in $(i)$ of Theorem~\ref{Bezdek-on-contact-numbers}. Namely, Conjecture~\ref{Bezdek-Molnar-lattice-estimate-II} implies in a straightforward way that
$$
{\rm svol}_2\left({\rm bd}\left(\bigcup_{i=1}^n(\mathbf{c}_i+2\mathbf{B})\right)\right)
$$
$$
\le 16\pi n - \frac{1}{6\left(1-\frac{\sqrt{3}}{2} \right)}16\pi \left(1-\frac{\sqrt{3}}{2} \right) C(n) =16\pi n - \frac{8\pi}{3}C(n)\ .
$$

The above inequality combined with Corollary~\ref{third-lower-bound-for-surface-area-in-3D} yields 
$$
4(18\pi)^{\frac{1}{3}}n^{\frac{2}{3}}< 16\pi n - \frac{8\pi}{3}C(n) \ ,
$$
from which the inequality  
$$
C(n)< 6n-\frac{3\sqrt[3]{18\pi}}{2\pi}n^{\frac{2}{3}}=6n-1.8326\dots n^{\frac{2}{3}}
$$
follows in a straightforward way.  This would improve the inequality $(i)$ of Theorem~\ref{Bezdek-on-contact-numbers} in a significant way.

\vspace{1cm}

\medskip

\noindent
K\'aroly Bezdek
\newline
Department of Mathematics and Statistics, University of Calgary, Canada,
\newline
Department of Mathematics, University of Pannonia, Veszpr\'em, Hungary,
\newline
and
\newline
Institute of Mathematics, E\"otv\"os University, Budapest, Hungary.
\newline
{\sf E-mail: bezdek@math.ucalgary.ca}

\end{document}